\documentclass[a4paper,11pt]{amsart}
\usepackage{paralist}
\usepackage{setspace}
\usepackage{geometry}
\usepackage{amsmath}
\usepackage{amssymb}
\usepackage{amsfonts}
\usepackage{color}
\usepackage{enumerate}
\usepackage{booktabs}
\usepackage{etoolbox}
\usepackage{verbatim}
\usepackage{graphicx}
\usepackage{subfig}
\usepackage{bm}
\usepackage{tikz}
\usepackage[bookmarks,bookmarksopen=false,                  pdfauthor={Autor},                  pdftitle={Titel},                  colorlinks=true,                  linkcolor=blue,                  citecolor=blue,                  urlcolor=blue,                ]{hyperref}
\usepackage{multicol}

\allowdisplaybreaks
\patchcmd{\thebibliography}{\section*{\refname}}{}{}{}
\bgroup

\newtheorem{theorem}{Theorem}

\newtheorem{proposition}[theorem]{Proposition}
\newtheorem{corollary}[theorem]{Corollary}

\newtheorem{lemma}[theorem]{Lemma}
\newtheorem{remark}[theorem]{Remark}

\renewcommand{\proof}{\textbf{Proof:\ }}
\newcommand{\pbox}{\hfill$\Box$\\}

\newcommand{\R}{\mathbb{R}}

\newcommand{\N}{\mathbb{N}}
\newcommand{\C}{\mathbb{C}}

\newcommand{\D}{\mathbb{D}}

\definecolor{darkviolet}{rgb}{0.58,0,0.83}

\usetikzlibrary{patterns}

\usepackage{currvita}
\begin{document}

\title[Large sieve principles for the wavelet
transform]{Donoho-Logan large sieve principles for the 
\\ wavelet transform }
\author[L. D. Abreu]{Lu\'{\i}s Daniel Abreu}
\address[L. D. A.]{Faculty of Mathematics \\
	University of Vienna \\
	Oskar-Morgenstern-Platz 1, A-1090 Vienna, Austria }
 \email{abreul22@univie.ac.at}

\author[M. Speckbacher]{Michael Speckbacher}

\address[M. S.]{Acoustics Research Institute\\ Austrian Academy of Sciences\\Dominikanerbastei 16, A-1010 Vienna, Austria }
\email{michael.speckbacher@oeaw.ac.at}
\date{}
\begin{abstract}
\noindent In this paper we formulate Donoho and Logan’s large sieve principle for the 
wavelet transform on the Hardy space, adapting the concept of maximum
Nyquist density to the hyperbolic geometry of the underlying space. The results
provide deterministic guarantees for $L_1$-minimization methods and hold for a class
of mother wavelets that constitutes an orthonormal basis of the Hardy space and
can be associated with higher hyperbolic Landau levels. Explicit calculations of the
basis functions reveal a connection with the Zernike polynomials. We prove a novel
local reproducing formula for the spaces in consideration and use it to derive concentration
estimates of the large sieve type for the corresponding wavelet transforms.
We conclude with a discussion of optimality of localization and Lieb inequalities in
the analytic case by building on recent results of Kulikov, Ramos and Tilli based on
the groundbreaking methods of Nicola and Tilli. This leads to a sharp uncertainty
principle and a local Lieb inequality for the wavelet transform.
\end{abstract}
\maketitle

\noindent \textbf{MSC2020:} 42C40, 46E15, 30H10, 42C15, 11N36\\
\noindent \textbf{Keywords:} {large sieve principle, wavelet transform, concentration estimates, maximum Nyquist rate, wavelet coorbit spaces, uncertainty principles}

\section{Introduction}

Let $d\mu ^{+}(z)= \text{Im}\left(z\right) ^{-2}d\mu _{\mathbb{C}^{+}}(z)$,
where $d\mu _{\mathbb{C}^{+}}(z)\ $is the Lebesgue measure in $\mathbb{C}^{+}
$. We denote the \emph{hyperbolic measure} of a  {nonempty} set $\Delta \subset \mathbb{C}^{+}$ by$\ |\Delta |_{h}:=\int_{\Delta }d\mu ^{+}(z)$ and define the weighted Lebesgue spaces on $\C^+$ via the norms $\|F\|_{L^p(\C^+)}^p=\int_{\C^+}|F(z)|^pd \mu^+(z)$. Let $0<R<1$, and $1\leq p<\infty$. We will be concerned with inequalities of the type 
\begin{equation}
\nu_p (\Delta )=\sup_{f\in \mathcal{C}o\, {L}^{p}(\mathbb{C}^{+})}\frac{%
\int_{\Delta }|W_{\psi }f(z)|^{p}d\mu ^{+}(z)}{\int_{\mathbb{C}^{+}}|W_{\psi
}f(z)|^{p}d\mu ^{+}(z)}\leq \frac{\rho (\Delta ,R)}{C_{\psi }(R)}\text{.}
\label{1}
\end{equation}
Here, $W_{\psi }f$ denotes the \emph{continuous analytic wavelet
transform} of a function $f$ with respect to a wavelet\ $\psi $, $\mathcal{C}%
o\, {L}^{p}(\mathbb{C}^{+})$ is the wavelet coorbit space with the
norm $\Vert f\Vert _{\mathcal{C}o\, {L}^{p}(\mathbb{C}^{+})}=\Vert
W_{\psi_0^\alpha}f\Vert _{ {L}^{p}(\mathbb{C}^{+})}$ (see the preliminaries
section for precise definitions) and $\rho (\Delta ,R)$ is the \emph{maximum
Nyquist density} given by 
\begin{equation}
\rho (\Delta ,R):=\sup_{z\in \mathbb{C}^{+}}|\Delta \cap \mathcal{D}%
_{R}(z)|_{h}\text{,}  \label{def-max-nyq}
\end{equation}%
where $\mathcal{D}_{R}(z)$ is a disc with respect to the pseudohyperbolic
distance on $\mathbb{C}^{+}$. Sets such that $\rho (\Delta ,R)\ll|\Delta |_{h}$%
\ will be said to be $R$-\emph{sparse in the hyperbolic measure}. 

The problem consists of obtaining the constant $C_{\psi }(R)$, which is
independent of $\Delta $ and grows with $R$. Given such a constant, the
estimate  \eqref{1} shows its full potential  for sets $\Delta $ which have low
concentration of measure in any hyperbolic disk $\mathcal{D}_{R}(z)$, allowing the use of $L_{1}$-minimization methods for signal recovery as in \cite{dosta89,dolo92,abspe18-sieve}. See also \cite{model,HrySpe,huli,JaSpe,nicola-lca} for similar concentration estimates.

The inequality \eqref{1} will be proven for an orthonormal basis $\{\psi
_{n}^{\alpha }\}_{n\geq 0}$ of the Hardy space $H^{2}(\mathbb{C}^{+})$
defined via the Fourier transform as 
\begin{equation*}
\widehat{\psi _{n}^{\alpha }}(t):={\sqrt{ \frac{ 2^{\alpha +2}\pi n!}{\Gamma
(n+\alpha +1)}}}\ t^{\frac{\alpha }{2}}e^{-t}L_{n}^{\alpha }(2t),\quad t>0,
\end{equation*}%
where $L_{n}^{\alpha }$, $ {\alpha>0,}$ denotes the generalized Laguerre polynomial of
degree $n$, 
\begin{equation*}
L_{n}^{\alpha }(t)=\sum_{k=0}^{n}\frac{(-1)^{k}}{k!}\binom{n+\alpha }{n-k}%
t^{k},\quad t>0.
\end{equation*}%
Note that the normalizing constant for $\psi _{n}^{\alpha }$ is chosen such
that  $\big\|\widehat{\psi
_{n}^{\alpha }}\big\|_{L^2(\R^+,t^{-1})}^2= {4\pi}/{\alpha }$, and $\Vert \psi _{n}^{\alpha }\Vert _{{^{2}}}^{2}=1$. This family was also the central object
of \cite{abspe22} where  conditions for frames generated by
orbits of Fuchsian groups in terms of a Nyquist-rate were obtained. The family $\{\psi _{n}^{\alpha }\}_{n\geq0}$ is connected to the eigenspaces of the Maass Laplacian (see Section 2.2), a Schr\"{o}dinger
operator that is of central importance in number theory due to its role in the
theory of Maass forms \cite{maass1949neue}, and Selberg trace formulas \cite{patterson1975laplacian}, as well as 
in physics, where the pure point spectrum of the Maass Laplacian admits an
interpretation as a hyperbolic analogue of the Euclidean Landau levels \cite{co87}. It
has been shown in \cite{abspe22} that these wavelet eigenspaces are $PSL(2,\mathbb{R}%
)$ invariant (a convenient property to define point processes \cite{affine}) and also that, assuming reasonably mild
restrictions on $\psi $ (see Theorem 3 in \cite{abspe22}), this is essentially the only choice with such property (see Theorem 3 in \cite{abspe22}). The sequence $\{\psi
_{n}^{\alpha }\}_{n\geq 0}$ also arises as the best localized functions of Daubechies and Paul's localization problem in the wavelet domain \cite{daupau88}, and have been used as time-scale tapers in spectral estimation \cite{bayram1996multiple} under the name of \emph{Morse wavelets}. It
should be noticed that the slightly different choice of wavelets whose Fourier transforms are a
constant multiple of $t^{\frac{\alpha +1}{2}}e^{-t}L_{n}^{\alpha }(2t)$ (note the different exponent of $t$)
leads to the polyanalytic Bergman  structure discovered by Vasilevski \cite{va99} (see \cite{hutnik2011wavelets,ab17,berge2022affine} for these special
choices in connection to poly-Bergman spaces, which are also fundamental in the theory of the affine Wigner distribution \cite{berge2022affineACHA} and in commutative algebras of Toeplitz operators \cite{vasilevski2008commutative}). 
What we would like to emphasize is that, in the STFT case, Hermite functions lead both to the polyanalytic structure and the Euclidean Landau levels \cite{abspe18-sieve}, but, for for the wavelet transform, \emph{the choices leading to Vasilevski's polyanalytic decomposition  structure \cite{va99} are different from those leading to the $PSL(2,\mathbb{R})$ invariant hyperbolic Landau level eigenspaces}. 

Our proof  of \eqref{1}  starts with a general Schur-type argument, as in \cite{abspe18-sieve}. Then the explicitly defined constant $C_{\psi
_{n}^{\alpha }}(R)=C_{n}^\alpha(R)$ is computed using the following local
reproducing formula for the wavelet coefficients: 
\begin{equation}
W_{\psi _{n}^{\alpha }}f(z)=\frac{1}{C_{n}^{\alpha }(R)}\int_{\mathcal{D}%
_{R}(z)}W_{\psi _{n}^{\alpha }}f(w)K_{\psi _{n}^{\alpha }}(z,w)d\mu
^{+}(w),\qquad z\in \mathbb{C}^{+}\text{.}  \label{local_rep}
\end{equation}%
Proving this formula will constitute a significant part of the paper. In \cite{abspe18-sieve}, a local reproducing formula for the short-time
Fourier transform (STFT) with Hermite windows was proven using the well-known correspondence between the  STFT and complex Hermite polynomials, which are orthogonal in
concentric discs of the plane. Since we could not find a similar
correspondence for wavelets in the literature, we\ have computed the wavelet
transforms $W_{\psi _{n}^{\alpha }}\psi _{m}^{\alpha }$, and the result
yields a correspondence (up to a conformal map) to the well known Zernike
polynomials in the disc \cite{wuensch05}, defined in terms of Jacobi polynomials in a fashion
reminiscent of the definition of the complex Hermite polynomials using
Laguerre functions \cite{ghanmi2008class,ismail2016analytic}. This   allows us to show the following orthogonality in
concentric pseudohyperbolic discs: 
\begin{equation*}
\int_{\mathcal{D}_{R}}W_{\psi _{n}^{\alpha }}\psi _{m}^{\alpha }(z)\overline{%
W_{\psi _{n}^{\alpha }}\psi _{k}^{\alpha }(z)}d\mu ^{+}(z)=C_{m,n}^{\alpha
}(R)\delta _{m-k}.
\end{equation*}%
The case $n=m$ of this orthogonality relation is then sufficient to show \eqref{local_rep}. This is a fundamental step in our derivation of the estimates of the maximal Nyquist rate, where the radial nature of the basis functions plays a role. Localization problems have been considered for other transforms (e.g., \cite{Mejjaoli2020Localization,mejjaoli2020spectral}), but it is not clear how to obtain local reproducing formulas in such settings.

For $p=1$, $\nu_1 (\Delta ) <\frac{1}{2}$ implies that%
\begin{equation}
\Vert W_{\psi }f\cdot \chi _{\Delta }\Vert _{{L}^{1}(\mathbb{C}%
^{+})}< \frac{1}{2}\Vert W_{\psi }f\Vert _{ {L}^{1}(%
\mathbb{C}^{+})}\text{,}  \label{sparse}
\end{equation}%
meaning that for \emph{every signal }$f\in \mathcal{C}o\, {L}^{1}(%
\mathbb{C}^{+})$, $W_{\psi }f(z)$\emph{\ is sparse }(low concentration on $%
\Delta $). By generalizing an observation of Donoho and Stark for
bandlimited discrete functions \cite[Theorem 9]{dosta89}, we obtain an
interesting reconstruction result. In the absence of noise, if $\nu_1(\Delta)<1/2$ and we only sense
the projection of a general $W_{\psi }f(z)$, $f\in \mathcal{C}o\, {L}%
^{1}(\mathbb{C}^{+})$ on $\Delta ^{c}$, then $W_{\psi }f(z)$ can be perfectly
reconstructed as the solution of the $L^{1}$-minimization problem%
\begin{equation*}
W_{\psi }f=\arg \min_{_{h\in \mathcal{C}o\, {L}^{1}(\mathbb{C}%
^{+})}}\Vert W_{\psi }h\Vert _{ {L}^{1}(\mathbb{C}^{+})}%
\text{, \ \ \ \ \ subject to } W_{\psi }h\big|_{\Delta^c}= W_{\psi }f\big|_{\Delta^c}%
\text{.}
\end{equation*}
Besides this result,  condition \eqref{sparse} allows several signal
approximation and recovery scenarios by using $L_{1}$-minimization \cite{abspe18-sieve}. 

A simple example shows that good estimates for $\nu_p (\Delta )$ can be obtained 
if $\Delta $ is sparse in the hyperbolic measure (low concentration in any disc $
\mathcal{D}_{R}(z)$). Consider first $
\Delta _{1}=\mathcal{D}_{R}$. Then $\rho (\Delta _{1},R)=|\mathcal{D}_{R}|_{h}
$. Now if $\Delta _{2}=\mathcal{D}_{R}\cup (\mathcal{D}_{R+2\delta }-\mathcal{D}_{R+\delta })$, with $0<\delta <R<1/2$, then $\rho (\Delta
_{2},R+\delta )=|\mathcal{D}_{R}|_{h}$. This gives the example of two sets $\Delta _{1},\Delta _{2}$ with $|\Delta
_{1}|_{h}<|\Delta _{2}|_{h}$ but with an estimate for $\nu_p (\Delta _{2})$
better than the estimate for $\nu_p 
 (\Delta _{1})$, since $\rho (\Delta
_{2},R+\delta )/C_{\psi }(R+\delta )<\rho (\Delta _{1},R)/C_{\psi }(R)$ (because, as we will see, in the estimates of the form \eqref{1}, $C_{\psi }(R)$ is independent of $\Delta $ and grows with $R$). 
This elementary observation already hints in the direction that motivated
Donoho and Logan, as well as our previous work \cite{abspe18-sieve}: if a
set is poorly concentrated in discs of big enough fixed radius $R$ (`$R$-\emph{sparse'}), then $\nu_p (\Delta )<\frac{1}{2}$.

Another interesting setting where the maximum Nyquist rate can be estimated is the hyperbolic Cantor set \cite{ABREU2022126699}. This has been done in 
 \cite[Section 3.1]{knutsen2022fractal} using our concept from \cite{abspe18-sieve}, and a similar calculation applies to the wavelet case. This approach may help in dealing with more general hyperbolic fractal sets, adapting the concept of porous sets from the time-frequency plane \cite{knutsen2022daubechies}.

The problem of estimating $\nu_p(\Delta)$ is closely connected to the study of  localization operators first introduced by Daubechies \cite{daub88} in the time-frequency domain and by Daubechies and Paul in the wavelet domain \cite{daupau88}. An analysis similar to the one in \cite{abspe18-sieve}
shows that our estimates fall short from being sharp. However, in the analytic  case $n=0$  it was recently shown by Ramos and Tilli \cite{RamosTilli} that
pseudohyperbolic disks maximize $\nu_2 (\Delta )$ among all sets $\Delta$ of a given hyperbolic measure, using methods from the
breakthrough paper \cite{Inventiones}, where the analogue result for the
time-frequency case has been obtained.
{For a discussion of optimal norm bounds for localization operators with non-binary symbols, see \cite{nicola2022norm,riccardi}.} We  include a section where we derive some direct consequences
of the  results in \cite{RamosTilli}, including a local Lieb inequality, an analogue of Theorem~5.2 in \cite{Inventiones}. 
This requires a Lieb-type
inequality for the Wavelet transform, that, as we shall see, has recently been proved in a slightly disguised
form by Kulikov \cite[Theorem~1.2]{kuli22}. 

The outline of the paper is as follows. In Section~\ref{sec:prel}, we gather all the concepts required to follow the paper in a relatively self-contained manner. 
In Section~\ref{sec:dbl-orth} we illustrate the  connection between the basis functions $W_{\psi _{n}^{\alpha }}\psi _{m}^{\alpha }$ and the Zernike polynomials. Then the double orthogonality of these basis functions in concentric pseudohyperbolic discs and the resulting local reproducing formula are proved. With these formulas at hand, we show our large sieve estimates  in Section~\ref{sec:LS}. Section~\ref{sec:opt} is devoted to optimality results in the analytic setting, 
concluding with a local Lieb uncertainty principle for analytic wavelets and a  short discussion about general Lieb inequalities. The computation of the basis functions $W_{\psi _{n}^{\alpha }}\psi _{m}^{\alpha }$, the explicit formulas for the reproducing kernels, and a proof of the integrability of the mother wavelets are left to the Appendix.

\section{Preliminaries}\label{sec:prel}

We use the following convention for the \emph{Fourier transform}: for $f\in
L^{1}(\mathbb{R})$, its Fourier transform $\widehat{f}$ is defined as 
\begin{equation*}
\widehat{f}(\xi )=\mathcal{F}(f)(\xi ):=\int_{\mathbb{R}}f(t)e^{-i\xi
t}dt,\qquad \xi \in \mathbb{R}.
\end{equation*}%
By standard arguments, the Fourier transform extends to $L^{2}(\mathbb{R})$
and \emph{Plancherel's formula} holds 
\begin{equation*}
\langle f_{1},f_{2}\rangle = {\frac{1}{2\pi}}\langle \widehat{f_{1}},\widehat{f_{2}}\rangle
,\qquad \text{ for }f_{1},f_{2}\in L^{2}(\mathbb{R}).
\end{equation*}%
We   use the basic notation for $H^{2}(\mathbb{C}^{+})$, the \emph{Hardy space on the upper half plane}, of analytic functions in $\mathbb{C}^{+}$ with the norm 
\begin{equation*}
\left\Vert f\right\Vert _{H^{2}(\mathbb{C}^{+})}=\text{ }\sup_{0<s<\infty }\int_{-\infty }^{\infty }\left\vert
f(x+is)\right\vert ^{2}dx<\infty \text{.}
\end{equation*}
To simplify the computations it is often convenient to use the equivalent
definition (since the Paley-Wiener theorem \cite{DGM} gives $\mathcal{F}
(H^{2}(\mathbb{C}^{+}))=L^{2}(\R^+ )$) 
\begin{equation*}
H^{2}(\mathbb{C}^{+})=\left\{ f\in L^{2}(\mathbb{R}):(\mathcal{F}f)(\xi )=0
\text{ for almost all }\xi <0\right\} \text{.}
\end{equation*}

\subsection{The Continuous Wavelet Transform }

\subsubsection{The affine group}

Consider the {\em $ax+b$ group} (see \cite[Chapter 10]{groe1} for the listed
properties) $G\sim \mathbb{R}\times \mathbb{R}^{+}\sim $ $\mathbb{C}^{+}$
with the multiplication
\begin{equation*}
(x,s)\cdot (x^{\prime },s^{\prime })=(x+sx^{\prime },ss^{\prime })\text{.}
\end{equation*}%
The identification $G\sim $ $\mathbb{C}^{+}$\ is done by setting $(x,s)\sim
x+is$. The neutral element of the group is $(0,1)\sim i$ and the inverse
element is given by $(x,s)^{-1}=(-\frac{x}{s},\frac{1}{s})\sim -\frac{x}{s}%
+\frac{i}{s}$. The $ax+b$ group is not\ unimodular, since the left Haar
measure on $G$ is $%
\frac{dxds}{s^2}$ and the right Haar measure on $G$ is $%
\frac{dxds}{s}$. \ Throughout this paper we will follow the convention that
for $z,w\in \mathbb{C}^{+}$, $zw$ denotes the product of two complex numbers
and $z\cdot w$ denotes the above group multiplication induced by the  affine group. 
For a set $\Delta \subset \mathbb{C}^{+}$ we define $\Delta ^{-1}:=\{z\in 
\mathbb{C}^{+}:\ z^{-1}\in \Delta \}$. The left Haar measure of a set $%
\Delta \subseteq G$ 
\begin{equation*}
|\Delta |=\int_{\Delta }\frac{dxds}{s^{2}}\text{, }
\end{equation*}
equals, under the identification of the $ax+b$ group with $\mathbb{C}^{+}
$,  the hyperbolic measure
\begin{equation*}
 |\Delta |=\left\vert \Delta \right\vert _{h}:= \int_{\Delta }
\text{Im}(z)^{-2}\,d\mu _{\mathbb{C}^{+}}(z)\text{,}
\end{equation*}
where $d\mu _{\mathbb{C}^{+}}(z)\ $is the Lebesgue measure in $\mathbb{C}^{+}
$. We will write 
\begin{equation}
d\mu ^{+}(z)=  \text{Im}\left(z\right) ^{-2}d\mu _{\mathbb{C}^{+}}(z)\text{.}
\label{measure}
\end{equation}

\subsubsection{The continuous wavelet transform}

For every $x\in \mathbb{R}$ and $s\in \mathbb{R}^{+}$, define the
translation operator  {$T_{x}:L^2(\R)\to L^2(\R)$} by $T_{x}f(t)=f(t-x)$, and the dilation operator   {$D_{s}:L^2(\R)\to L^2(\R)$} by $D_{s}f(t)=\frac{1
}{\sqrt{s}}f(t/s)$. Let $z=x+is\in  \mathbb{C}^{+}$ and define the unitary representation $\pi:G\to \mathcal{U}(H^2(\C^+))$
\begin{equation}
\pi (z)\psi (t):=T_{x}D_{s}\psi (t)=s^{-\frac{1}{2}}\psi (s^{-1}(t-x)),\qquad \psi \in H^{2}(\mathbb{C}^{+}) .\label{representation}
\end{equation}
To properly define the wavelet transform, one needs to take into account that the $ax+b$ group is
not\ unimodular
which requires the following additional condition on the integrability of 
$\psi\in H^2(\C^+)$ 
\begin{equation}
0<\left\Vert \mathcal{F}\psi \right\Vert _{L^{2}(\mathbb{R}^{+},t^{-1}dt)}^{2}=:C_{\psi }<\infty. \label{Adm_const}
\end{equation}
Functions satisfying \eqref{Adm_const} are called \emph{admissible }and the
constant $C_{\psi }$ is the \emph{admissibility constant}.  
The \emph{continuous analytic wavelet transform} of a function $f$ with
respect to a wavelet\ $\psi $ is defined as 
\begin{equation}
W_{\psi }f(z)=\left\langle f,\pi ({z})\psi \right\rangle _{H^{2}\left( 
\mathbb{C}^{+}\right) }=\frac{1}{\sqrt{s}}\int_{\R}\tilde{f}(t)\overline{\tilde{\psi} \left( \frac{t-x}{s}\right)}dt,\quad z=x+is\in \mathbb{C}^{+}, 
\end{equation}
where $\tilde{f}(t)=\lim_{y\to0^+}f(t+iy)$.
Using $\mathcal{F}(H^{2}(\mathbb{C}^{+}))=L^{2}(\R^+ )$, this can also
be written (and we will do it as a rule to simplify the calculations) as%
\begin{equation}
W_{\psi }f(z)={\frac{\sqrt{s}}{2\pi}}\int_{\R^+ }(\mathcal{F}f)(\xi )e^{ix\xi } \overline{(\mathcal{F}%
\psi)}(s\xi )d\xi \text{.}  \label{anwav}
\end{equation}
As proven recently in \cite{AnalyticWavelet}, $W_{\psi }f(z)$ only leads to
analytic (Bergman) phase spaces for a very special choice of $\psi $, but it
is common practice to also call it  continuous analytic wavelet
transform for general wavelets $\psi$. The orthogonality relation
\begin{equation}
\int_{\mathbb{C}^{+}}W_{\psi _{1}}f_{1}(z)\overline{W_{\psi _{2}}f_{2}(z)
}d\mu ^{+}(z)= \left\langle \mathcal{F}\psi _{1},\mathcal{F}\psi
_{2}\right\rangle _{L^{2}(\mathbb{R}^{+},t^{-1})}\left\langle f_{1},f_{2}\right\rangle _{H^{2}\left( \mathbb{C}
^{+}\right) }\text{,}  \label{ortogonalityrelations}
\end{equation}
is valid for all $f_{1},f_{2}\in H^{2}\left( \mathbb{C}^{+}\right) $ and $
\psi _{1},\psi _{2}\in H^{2}\left( \mathbb{C}^{+}\right) $ admissible. Setting $\psi _{1}=\psi _{2}=\psi $\ and $f_{1}=f_{2}$\ in \eqref{ortogonalityrelations} then gives
\begin{equation}
\int_{\mathbb{C}^{+}}\left\vert W_{\psi }f(z)\right\vert ^{2}d\mu
^{+}(z)=C_{\psi }\left\Vert f\right\Vert _{H^{2}\left( \mathbb{C}^{+}\right)
}^{2},  \label{isometry}
\end{equation}
that is, the continuous wavelet transform is a multiple of an isometric inclusion $
W_{\psi }:H^{2}\left( \mathbb{C}^{+}\right) \rightarrow L^{2}(\mathbb{C}^{+}
)$. Setting $\psi _{1}=\psi _{2}=\psi $\ and $f_{2}=\pi
(z)\psi $ in \eqref{ortogonalityrelations} also shows that for every $f\in H^{2}\left( 
\mathbb{C}^{+}\right) $ one has 
\begin{equation}
W_{\psi }f(z)=\frac{1}{C_{\psi }}\int_{\mathbb{C}^{+}}W_{\psi }f(w)\langle
\pi (w)\psi ,\pi (z)\psi \rangle d\mu ^{+}(z),\quad z\in \mathbb{C}^{+}, \label{eq:rep-eq}
\end{equation}
where the integral converges in the weak sense.
Thus, the range of the wavelet transform 
\begin{equation*}
W_{\psi }\left( H^{2}\left( \mathbb{C}^{+}\right) \right) :=\big\{F\in L^{2}(
\mathbb{C}^{+},\mu ^{+}):\ F=W_{\psi }f,\ f\in H^{2}\left( \mathbb{C}
^{+}\right) \big\}
\end{equation*}
is a reproducing kernel subspace of $L^{2}(\mathbb{C}^{+} )$ with
kernel
\begin{equation}
K_{\psi }(z,w)=\frac{1}{C_{\psi }}\langle \pi (w)\psi ,\pi (z)\psi \rangle
_{H^{2}\left( \mathbb{C}^{+}\right) }=\frac{1}{C_{\psi }}W_{\psi }\psi ({w}
^{-1}\cdot z),  \label{RepKern}
\end{equation}
and $K_{\psi }(z,z)= {\Vert \psi \Vert _{2}^{2}
}/{C_{\psi }}$\text{.}
The Fourier transform $\mathcal{F}:H^{2}\left( \mathbb{C}^{+}\right)
\rightarrow L^{2}(\R^+ )$ can be used to simplify computations, since for $z=x+is$ and $w=x'+is'$
\begin{align}
\left\langle \pi (w)\psi ,\pi (z)\psi \right\rangle _{H^{2}\left( \mathbb{C}%
^{+}\right) }&= {\frac{1}{2\pi}}\left\langle \widehat{\pi (w)\psi },\widehat{\pi (z)\psi }%
\right\rangle _{L^{2}(\mathbb{R}^{+})}\notag
\\
&= {\frac{1}{2\pi}}\left( ss^{\prime }\right) ^{\frac{%
1}{2}}\int_{\R^+ }\widehat{\psi }(s^{\prime }\xi )\overline{\widehat{%
\psi }\left( s\xi \right) }e^{i(x-x^{\prime })\xi }d\xi \text{.}
\label{kernelF}
\end{align}

\subsubsection{Wavelet coorbit spaces}

It is commonly known that the representation $\pi $ is also \emph{integrable}, that
is, there exist admissible mother wavelets $\psi \neq 0$ such that 
\begin{equation}
\int_{\mathbb{C}^{+}}|\langle \psi ,\pi (z)\psi \rangle |d\mu ^{+}(z)<\infty
.  \label{eq:integr}
\end{equation}
We will show in Appendix~\ref{app:int} that we can choose $\psi=\psi_n^\alpha$, $\alpha>1$, here. 
For such a mother wavelet, we define the space of test functions 
\begin{equation*}
\mathcal{H}_{1}:=\big\{ f\in H^{2}\left( \mathbb{C}^{+}\right) :\ W_{\psi}f\in
L^{1}(\C^+)\big\} ,
\end{equation*}%
and denote by $\mathcal{H}_{1}^{\prime }$ its anti-dual space, i.e., the
space of antilinear continuous functionals on $\mathcal{H}_{1}$. Then the 
\emph{coorbit space} with respect to ${L}^{p}(\mathbb{C}^{+})$, $%
1\leq p\leq \infty $, is defined as 
\begin{equation*}
\mathcal{C}o\,{L}^{p}(\mathbb{C}^{+}):=\left\{ f\in \mathcal{H}%
_{1}^{\prime }:W_{\psi}f\in  {L}^{p}(\mathbb{C}^{+})\right\} ,
\end{equation*}%
equipped with the natural norm $\Vert f\Vert _{\mathcal{C}o\,{L}^{p}(%
\mathbb{C}^{+})}^p:=\int_{\C^+}| W_{\psi}f(z)|^pd \mu^+(z)$, $1\leq p<\infty$, and the usual modification if $p=\infty$.
See, e.g., \cite{fegr89,fegr89II,groe91} for more information on coorbit
space theory. The spaces $\mathcal{C}o\, {L}^{p}(\mathbb{C}^{+})$ are
Banach spaces for $1\leq p\leq \infty $, they are independent of the
particular choice of mother wavelet $\psi\in \mathcal{H}_{1}$, the reproducing
formula \eqref{eq:rep-eq} extends to $\mathcal{C}o\, {L}^{p}(%
\mathbb{C}^{+})$, and 
\begin{equation*}
\mathcal{C}o\, {L}^{1}(\mathbb{C}^{+})=\mathcal{H}_{1}\subset 
\mathcal{C}o\, {L}^{2}(\mathbb{C}^{+})=H^{2}\left( \mathbb{C}%
^{+}\right) \subset \mathcal{C}o\, {L}^{\infty }(\mathbb{C}^{+})=%
\mathcal{H}_{1}^{\prime }.
\end{equation*}

\subsection{Hyperbolic Landau Level Spaces}

\label{sec:Landau-levels}

We first make an important remark on the transforms $W_{\psi _{n}^{\alpha
}}$. \ The spaces $W_{\psi _{n}^{\alpha }}(H^2(\C^+))$, $n\geq 0$, are not orthogonal for a
general choice of  $\alpha $ because the wavelets $\psi _{n}^{\alpha }$ are orthogonal in
the Hardy space and not in the space of admissible wavelets. But if we make
the particular choice $\alpha =2B-2n-1$, they become orthogonal, that is, for $%
n,m\in \big\{0,1,\ldots ,\lfloor B-\frac{1}{2}\rfloor \big\}$, 
\begin{equation*}
\langle W_{\psi _{n}^{2B-2n-1}}f,W_{\psi _{m}^{2B-2m-1}}g\rangle _{{L%
}^{2}(\mathbb{C}^{+})}=\frac{4\pi}{2B-2n-1}  {\langle f,g\rangle_{H^2(\C^+)}}\delta _{n,m},\quad {f,g\in H^2(\C^+)}.
\end{equation*}%
This is related to the identification of the spaces $W_{\psi _{n}^{2B-2n-1}}(H^{2}(\mathbb{C}^{+}))$ as the eigenspaces of the Maass-Landau levels
operator with a constant magnetic field $B$ originally studied in number theory in the theory of Maass forms \cite{maass1949neue}, and the Selberg trace formula \cite{patterson1975laplacian}. We provide a brief account of this Schr\"{o}dinger operator from a physical perspective.

The \emph{Schr{\"{o}}dinger operator} describing the dynamics of a charged
particle moving in $\mathbb{C}^{+}$ under the action of the constant magnetic field $B
$ \cite{co87} is given by 
\begin{equation}
H_{B}=-s^{2}\left( \frac{\partial ^{2}}{\partial x^{2}}+\frac{\partial ^{2}}{%
\partial s^{2}}\right) +2iBs\frac{\partial }{\partial x}\text{.}
\label{eq:schroedinger-op}
\end{equation}%
$H_{B}$ is an elliptic, and densely defined operator on ${L}^{2}(%
\mathbb{C}^{+})$. Its spectrum consists of a continuous part $[1/4,\infty )$
corresponding to \emph{scattering states} and a finite number of eigenvalues 
\begin{equation*}
\lambda _{n}^{B}=(B-n)(B-n-1),\qquad n=0,1\ldots ,\Big\lfloor B-\frac{1}{2}
\Big\rfloor \text{.}
\end{equation*}
The eigenvalues exist provided that $2B>1$, and each eigenvalue corresponds
to an infinite dimensional reproducing kernel Hilbert space $E_{n}^{B}(
\mathbb{C}^{+})$. From a physical viewpoint, this condition guarantees
that the magnetic field is strong enough to capture particles in a closed
orbit, see \cite[p. 189]{gro88}. The eigenspaces corresponding to $\lambda
_{n}^{B}$ are called \emph{hyperbolic Landau levels} and we define them by 
\begin{equation*}
E_{n}^{B}(\mathbb{C}^{+}):=\big\{F\in  {L}^{2}(\mathbb{C}^{+}):\
H_{B}F=\lambda _{n}^{B}F\big\}\text{.}
\end{equation*}
In \cite[Section 4,5]{affine} it is shown in detail (this connection was first observed by Mouayn \cite{mouayn2003characterization}) that 
\begin{equation*}
W_{\psi _{n}^{2B-2n-1}}(H^{2}(\mathbb{C}^{+}))=E_{n}^{B}(\mathbb{C}^{+})
\text{.}
\end{equation*}
This way we obtain an orthogonal basis for all the spaces $E_{n}^{B}(\mathbb{C}
^{+})$, by noting that for $n,m\in
\{0,1,\ldots ,\lfloor B-\frac{1}{2}\rfloor \}$, \eqref{ortogonalityrelations} gives 
\begin{equation*}
\langle W_{\psi _{n}^{2B-2n-1}}\psi _{k}^{\alpha },W_{\psi
_{m}^{2B-2m-1}}\psi _{l}^{\alpha }\rangle _{{L}^{2}(\mathbb{C}^{+})}=%
\frac{2}{2B-2n-1}\delta _{n,m}\delta _{k,l}\text{.}
\end{equation*}
Our large sieve inequalities for this case provide conditions allowing for the recovery of a continuous coherent state in higher hyperbolic Landau levels, from partial information, using $L_{1}$-minimization, as outlined in the introduction.
\subsection{Pseudohyperbolic Metric and M{\"o}bius Transformation}

The \emph{pseudohyperbolic distance} on $\mathbb{C}^{+}$ is given by 
\begin{equation*}
\varrho (z,w):=\left\vert \frac{z-w}{z-\overline{w}}\right\vert ,\qquad
z,w\in \mathbb{C}^{+},
\end{equation*}%
and the \emph{pseudohyperbolic disk} of radius $R>0$ centered at $z\in 
\mathbb{C}^{+}$ is defined by $\mathcal{D}_{R}(z):=\{\omega \in \mathbb{C}%
^{+}:\ \varrho (z,w)<R\}$. We write for short $\mathcal{D}_{R}=\mathcal{D}%
_{R}(i)$. Note that $\varrho $ only takes values in the half open interval $
[0,1)$. It can be checked that the pseudohyperbolic distance is symmetric
and  invariant under the action of the affine group. This leads to several
useful properties, like $\varrho (z,w)=\varrho (z^{-1}\cdot w,i)$ and $
\varrho (u,z\cdot w)=\varrho (z^{-1}\cdot u,w)$.

Let $\mathbb{D}_{R}$ denote the disk of radius $R>0$ in $\mathbb{C}$. We
write for short $\mathbb{D}:=\mathbb{D}_{1}$ to denote the unit disk. The
mapping $T:\mathbb{D}\rightarrow \mathbb{C}^{+}$,  
\begin{equation*}
T(u):=i\frac{1+u}{1-u},\quad u\in \mathbb{D}
\end{equation*}%
maps the pseudohyperbolic distance in $\mathbb{C}^{+}$ to the
pseudohyperbolic distance in $\mathbb{D}$:
\begin{equation*}
\varrho (T(z),T(w))=\left\vert \frac{w-z}{1-\overline{z}w}\right\vert
=\varrho _{\mathbb{D}}(z,w)\text{.}
\end{equation*}
The mapping $T$ will be used to map some integrals from the upper half plane to
the unit disk according to the change of variables formula
\begin{equation}\label{eq:change-of-variable}
\int_{\mathbb{C}^{+}}F(z)s^{\beta }d\mu_{\C^+}(z)=\int_{\mathbb{D}}F(T(u))\frac{
4(1-|u|^{2})^{\beta }}{|1-u|^{2\beta +4}}du,\qquad \beta \in \mathbb{R},\ z=x+is.
\end{equation}

\section{Double Orthogonality and the Local Reproducing Formula}\label{sec:dbl-orth}

The following family of   polynomials defined on $[0,1]$ is closely related to  the \emph{ Zernike
polynomials}  \cite{wuensch05}, a family of orthogononal polynomials in $(z,\overline{z})$ defined on $\D$,
\begin{equation}
Z_{n,m}^{\alpha }(t):=\sqrt{\dfrac{\Gamma (\max (n,m)+\alpha
+1)\min (n,m)!}{\Gamma (\min (n,m)+\alpha +1)\max (n,m)!}}(-t)^{-\min
(n,m)}P_{\min (n,m)}^{(|n-m|,\alpha )}\left( 1-2t\right),
\label{eq:def-Phi_nm}
\end{equation}
where $P_n^{(\alpha,\beta)}$ denotes the \emph{Jacobi polynomials} which are given explicitly in \eqref{eq:jacobi}.
We will see in this section that 
\begin{equation*}
Z_{n,m}^{\alpha }\left( \left\vert \dfrac{z-i}{z+i}\right\vert ^{2}\right),
\quad z\in \mathbb{C}^{+},
\end{equation*}
plays a role in wavelet analysis similar to the one of complex Hermite
polynomials in time-frequency analysis \cite{abspe18-sieve}. The following result is proved in the
Appendix.

\begin{proposition}
\label{lem:explicit-wavelet-trafo} Let $\alpha>0$, and $z=x+is$ and $w=x^{\prime
}+is^{\prime }$ be in $\mathbb{C}^{+}$. For every $n\in \mathbb{N}_{0}$,
one has   {
\begin{equation}\label{C_psi_n^alpha}
C_{\psi_n^\alpha}=\frac{4\pi}{\alpha},
\end{equation}
}
and
\begin{equation}\label{eq:kernel-expl}
K_{\psi _{n}^{\alpha }}(z,w)=\frac{\alpha}{4\pi} \left( \frac{w-\overline{z}}{z-
\overline{w}}\right) ^{n}\left( \frac{{2}\sqrt{ss^{\prime }}}{i(\overline{w}-z)}
\right) ^{\alpha +1}P_{n}^{(0,\alpha )}\left( 1-2\left\vert \frac{z-w}{z-
\overline{w}}\right\vert ^{2}\right).
\end{equation}
 {Moreover, for $n,m\in\N_0$, it holds}
\begin{equation}
W_{\psi _{n}^{\alpha }}\psi _{m}^{\alpha }(z)= \left(\frac{z-i}{z+i}\right)^{m}\left(\frac{\overline{z}%
+i}{z+i}\right)^{n}\left( \dfrac{{2}\sqrt{s}}{1-iz}\right) ^{\alpha
+1}Z_{n,m}^{\alpha }\left( \left\vert \dfrac{z-i}{z+i}\right\vert
^{2}\right) .  \label{eq:explicit-wavelet}
\end{equation}
\end{proposition}

\begin{remark}
 For $n=0$, i.e. the case of the Cauchy wavelet, one has that the
Jacobi polynomial $P_0^{(0,\alpha)}$ is constantly  one. Consequently, 
\begin{equation}
\bm{\Phi}_{m}^{\alpha }(z):=s^{-\frac{\alpha +1}{2}}W_{\psi _{0}^{\alpha
}}\psi _{m}^{\alpha }(z)= \sqrt{\frac{\Gamma (m+\alpha +1)}{
\Gamma (\alpha +1)m!}}\left( \frac{z-i}{z+i}\right) ^{m}\left( \frac{{2}}{1-iz}
\right) ^{\alpha +1},\label{eq:n=0}
\end{equation}
which are (up to a constant factor due to different normalizations) the
basis functions of the Bergman space $
A_{\alpha -1}^2(\mathbb{C}^{+})$, where 
$$
A_{\alpha }^p(\mathbb{C}^{+}):=\{F\text{ holomorphic}:\ \|F\|_{A^p_\alpha(\C^+)}<\infty\},
$$
and 
\begin{equation}\label{bergman:halfplane}
\|F\|_{A_{\alpha }^p(\mathbb{C}^{+})}^p:=\int_{\C^+}|F(z)|^p s^{\alpha}d \mu_{\C^+}(z),
\end{equation}  
see,
e.g., \cite[Section 4.4]{abdoe12}.

\end{remark}

In the following, we show a local orthogonality relation for the family $%
\{W_{\psi _{n}^{\alpha }}\psi _{m}^{\alpha }\}_{m\in \mathbb{N}_{0}}$ and
subsequently derive a local reproducing formula.

\begin{theorem}
\label{thm:local-orthogonal} Let $\alpha >0$, $0<R<1$, and $n,m,k\in \mathbb{
N}_{0}$. Then it holds 
\begin{equation}
\int_{\mathcal{D}_{R}}W_{\psi _{n}^{\alpha }}\psi _{m}^{\alpha }(z)\overline{
W_{\psi _{n}^{\alpha }}\psi _{k}^{\alpha }(z)}d\mu ^{+}(z)=C_{m,n}^{\alpha
}(R)\delta _{m-k},  \label{eq:local-orthogonal}
\end{equation}
with 
\begin{equation*}
C_{n,m}^{\alpha }(R)=   {4\pi} \int_{0}^{R^{2}}r^{n+m}(1-r)^{\alpha
-1}Z_{n,m}^{\alpha }(r)^{2}dr.
\end{equation*}
If $n=m$, we write $C_{n}^{\alpha }(R)=C_{n,n}^{\alpha }(R)$ which is given
explicitly by 
\begin{equation}
C_{n}^{\alpha }(R)=  {4\pi}\int_{0}^{R^{2}}(1-r)^{\alpha
-1}P_{n}^{(0,\alpha )}(1-2r)^{2}dr.  \label{eq:def_Cn}
\end{equation}
\end{theorem}

\noindent \textbf{Proof:\ }
First we note that putting $z=T(u)$ in \eqref{eq:explicit-wavelet} yields
\begin{equation}\label{eq:aux-moebius}
W_{\psi _{n}^{\alpha }}\psi _{m}^{\alpha }(T(u))={(-1)^{n+\frac{\alpha+1}{2}}}u^{m}\overline{u}^{n}\left( \dfrac{1-u}{|1-{u}|}\right) ^{2n+\alpha +1}\left(  {1-|u|^{2}} \right) ^{\frac{\alpha +1}{2}}Z _{n,m}^\alpha(|u|^{2}).
\end{equation}
Using \eqref{eq:change-of-variable} for $\beta=-2$ then leads us to
\begin{align*}
\int_{\mathcal{D}_{R}}W_{\psi _{n}^{\alpha }}\psi _{m}^{\alpha }(z)\overline{%
W_{\psi _{n}^{\alpha }}\psi _{k}^{\alpha }(z)}d\mu ^{+}(z)& =\int_{\mathbb{D}%
_{R}}W_{\psi _{n}^{\alpha }}\psi _{m}^{\alpha }(T(u))\overline{W_{\psi
_{n}^{\alpha }}\psi _{k}^{\alpha }(T(u))}\frac{4du}{(1-|u|^{2})^{2}} \\
&= {4}\int_{\mathbb{D}_{R}}u^{m}\overline{u}^{k}|u|^{2n}(1-|u|^{2})^{%
\alpha -1}Z_{n,m}^{\alpha }(|u|^{2})Z_{n,k}^{\alpha }(|u|^{2})du \\
& =\circledast .
\end{align*}%
Integration using polar coordinates then yields 
\begin{align*}
\circledast & = {4}\int_{0}^{2\pi }\int_{0}^{R}e^{i(m-k)\varphi
}r^{m+k+2n+1}(1-r^{2})^{\alpha -1}Z_{n,m}^{\alpha }(r^{2})Z_{n,k}^{\alpha
}(r^{2})drd\varphi \\
& = {8}\pi \delta _{m-k}\int_{0}^{R}r^{2(n+m)+1}(1-r^{2})^{\alpha
-1}Z_{n,m}^{\alpha }(r^{2})^{2}dr \\
& ={4}\pi \delta _{m-k}\int_{0}^{R^{2}}r^{n+m}(1-r)^{\alpha -1}Z_{n,m}^{\alpha
}(r)^{2}dr,
\end{align*}%
which concludes the proof once we recall the definition of $Z_{n,m}^{\alpha }
$ in \eqref{eq:def-Phi_nm}. \hfill $\Box $\newline

Let us shortly consider the case  $n=0:$%
\begin{align*}
C_{0}^{\alpha }(R) &= {4\pi}\int_{0}^{R^{2}}(1-r)^{\alpha
-1}P_{0}^{(0,\alpha )}(1-2r)^{2}dr= {4\pi}\int_{0}^{R^{2}}(1-r)^{\alpha -1}dr \\
&= {\frac{4\pi}{\alpha}}\left( 1-(1-R^{2})^{\alpha }\right),
\end{align*}%
which converges to  {$C_{\psi_0^\alpha}=4\pi/\alpha$} (as $R\to 1$) as one expects from 
\eqref{ortogonalityrelations}.

\begin{theorem}
\label{thm:loc-rep} Let $n\in \mathbb{N}_{0}$, $\alpha >1$, and $0<R<1$. The
following identity holds in the weak sense 
\begin{equation}
\int_{\mathcal{D}_{R}}W_{\psi _{n}^{\alpha }}\psi _{n}^{\alpha }(z)\pi
(z)\psi _{n}^{\alpha }d\mu ^{+}(z)=C_{n}^{\alpha }(R)\ \psi _{n}^{\alpha }.
\label{eq:weak-reproducing}
\end{equation}
If $f\in \mathcal{C}o\,{L}^{\infty }(\mathbb{C}^{+})$, then the wavelet
coefficients can locally be reconstructed via 
\begin{equation}
W_{\psi _{n}^{\alpha }}f(z)=\frac{4\pi}{\alpha C_{n}^{\alpha }(R)}\int_{\mathcal{D}
_{R}(z)}W_{\psi _{n}^{\alpha }}f(w)K_{\psi _{n}^{\alpha }}(z,w) d\mu
^{+}(w),\qquad z\in \mathbb{C}^{+}.  \label{eq:strong-reprod}
\end{equation}
\end{theorem}

\noindent \textbf{Proof:\ } Setting $m=n$ in \eqref{eq:local-orthogonal} shows that 
\eqref{eq:weak-reproducing} holds weakly. If we apply $\pi (w)$ on both
sides of \eqref{eq:weak-reproducing} and subsequently take the inner product
with $f$, we obtain 
\begin{align*}
W_{\psi _{n}^{\alpha }}f(z)& =\langle f,\pi (z)\psi _{n}^{\alpha }\rangle \\
& =\frac{1}{C_{n}^{\alpha }(R)}\int_{\mathcal{D}_{R}}\langle f,\pi (z\cdot
w)\psi _{n}^{\alpha }\rangle \langle \pi (w)\psi _{n}^{\alpha },\psi
_{n}^{\alpha }\rangle d\mu ^{+}(w) \\
& =\frac{1}{C_{n}^{\alpha }(R)}\int_{\mathcal{D}_{R}(z)}\langle f,\pi
(w)\psi _{n}^{\alpha }\rangle \langle \pi (z^{-1}\cdot w)\psi _{n}^{\alpha
},\psi _{n}^{\alpha }\rangle {d\mu^+(w)} \\
& =\frac{1}{C_{n}^{\alpha }(R)}\int_{\mathcal{D}_{R}(z)}\langle f,\pi
(w)\psi _{n}^{\alpha }\rangle \langle \pi (w)\psi _{n}^{\alpha },\pi (z)\psi
_{n}^{\alpha }\rangle d\mu ^{+}(w) \\
& =\frac{ {C_{\psi_n^\alpha}}}{  C_{n}^{\alpha }(R)}\int_{\mathcal{D}_{R}(z)}W_{\psi _{n}^{\alpha
}}f(w)K_{\psi _{n}^{\alpha }}(z,w)d\mu ^{+}(w),
\end{align*}
where we used  $
\varrho (z,z\cdot w)=\varrho (i,w)$, and the left
invariance of the Haar measure. The second equality holds since the integral
and the duality pairing may be interchanged knowing that $f\in \mathcal{C}o\,
 {L}^{\infty }(\mathbb{C}^{+})$, and $K_{\psi _{n}^{\alpha }}(z,\, \cdot\, )\in 
 {L}^{1}(\mathbb{C}^{+})$.   {The result follows once we recall that $C_{\psi_n^\alpha}=\frac{4\pi}{\alpha}$.} \hfill $\Box $\newline

\section{Large Sieve Estimates}\label{sec:LS}

\subsection{A Schur-type estimate}

Let $(X,\mu )$, $(X,\nu )$ be measure spaces and $B_{1}\subset L^{1}(X,\nu )$
a Banach space. In \cite[Proposition 1]{abspe18-sieve}, a bound
on the embedding $(B_{1},\Vert \cdot \Vert _{L_{\nu }^{1}})\hookrightarrow
(B_{1},\Vert \cdot \Vert _{L_{\mu }^{1}})$ was derived using an argument
similar to Schur's test which we shortly recall here.

\begin{lemma}
\label{prop-measure-bound} Let $\mu $ be a positive $\sigma $-finite measure
on $X$, $B_{1}\subset L^{1}(X,\nu)$ be a Banach space, and $K:X\times X\rightarrow \mathbb{C} $
be 
such that $\mathcal{K}:B_{1}\rightarrow B_{1}$, 
\begin{equation*}
\mathcal{K}F(x):=\int_{X}F(y)K(x,y)d\nu(y),
\end{equation*}
is bounded and boundedly invertible on $B_{1}$. Then, for every $F\in B_{1}$, we have 
\begin{equation} 
\frac{\int_{X}|F|d\mu }{\int_X |F|d\nu}\leq \theta (K) \, \sup_{y\in
X}\int_{X}|K(x,y)|d\mu (x),  \label{eq-norm-bound2}
\end{equation}
where 
\begin{equation*}
\theta  (K):=\sup_{H \in B_{1}}\left( \dfrac{\Vert H \Vert _{L^1_\nu}}{\Vert \mathcal{K}H \Vert _{L^1_\nu}}\right) .
\end{equation*}
\end{lemma}

\subsection{Estimates with explicit constants}

Recall that the \emph{maximum Nyquist density} is given by 
\begin{equation*}
\rho (\Delta ,R):=\sup_{z\in \mathbb{C}^{+}}|\Delta \cap \mathcal{D}_{R}(z))|_{ {h}},
\end{equation*}
where $|\Delta |_{ {h}}:= \int_{\Delta }d\mu ^{+}(z)$ denotes \emph{hyperbolic
measure} of $\Delta $. We will also use a second notion of density given
by 
\begin{equation}
D_{n}^{\alpha }(\Delta ,R):= \sup_{z\in \mathbb{C}^{+}}\int_{\Delta \cap 
\mathcal{D}_{R}(z)}|\langle \pi(w)\psi_{n}^{\alpha },\pi(z)\psi_n^\alpha\rangle| {d\mu^+(w)}.
\label{Ar-density}
\end{equation}
Note that $D_{n}^{\alpha }(\Delta ,R)\leq \rho (\Delta ,R)$, as 
\begin{equation}
|\langle \pi(w)\psi_{n}^{\alpha },\pi(z)\psi_n^\alpha\rangle|\leq \Vert \psi _{n}^{\alpha }\Vert ^{2}_{H^2(\C^+)}=1.
\label{KernelEst}
\end{equation}
We have gathered all the ingredients to formulate and prove our main result.

\begin{theorem}
\label{main} Let $\Delta \subset \mathbb{C}^{+}$, $f\in \mathcal{C}o\,
 {L}^{p}(\mathbb{C}^{+})$, $1\leq p<\infty $, and   {$\alpha>1$.} For every $0<R<1$, it
holds 
\begin{equation}
\dfrac{\big\| W_{\psi _{n}^{\alpha }}f\big|_{\Delta }\big\|_{L^p(\C^+)}^p}{\Vert W_{\psi _{n}^{\alpha }}f\Vert _{L^p(\C^+)}^{p}}\leq  \frac{D_{n}^{\alpha
}(\Delta ,R)}{C_{n}^{\alpha }(R)}\leq   \frac{\rho (\Delta ,R)}{C_{n}^{\alpha
}(R)}\text{.}\   \label{theo1-eq}
\end{equation}
\end{theorem}

\noindent \textbf{Proof:\ } We take $K(z,w):=\langle \pi(w)\psi_{n}^{\alpha },\pi(z)\psi_n^\alpha\rangle \chi _{\mathcal{D}_{R}(z)}(w)$ and $B_1=W_{\psi_n^\alpha}\big(Co L^1(\C^+)\big)$ in Lemma~\ref{prop-measure-bound}. Then by Theorem~\ref{thm:loc-rep} we get
\begin{equation*}
\theta(K)=\sup_{f\in \mathcal{C}o\, {L}^{1}(\mathbb{C}^{+} )}\left( \dfrac{\Vert W_{\psi _{n}^{\alpha }}f\Vert _{L^{1}(\C^+ )}}{\big\|
\int_{\mathcal{D}_{R}(\cdot)}W_{\psi _{n}^{\alpha }}f(w)\langle \pi(w)\psi_{n}^{\alpha },\pi(\cdot)\psi_n^\alpha\rangle d\mu ^{+}(w)\big\| _{L^{1}(\C^+)}}\right) =\frac{1}{C_{n}^{\alpha }(R)}.
\end{equation*}
Thus, if $d\mu (z)=\chi _{\Delta }(z)d\mu^+(z)$, we have by \eqref{eq-norm-bound2} and \eqref{KernelEst}
\begin{equation*}
\dfrac{\big\| W_{\psi _{n}^{\alpha }}f\big|_{\Delta }\big\| _{ {L}^{1}(\C^+)}}{\Vert W_{\psi }^{\alpha }f\Vert _{ {L}^{1}(\C^+)}}\leq \frac{1}{C_{n}^{\alpha }(R)}\, \sup_{z\in \mathbb{C}^{+}}\int_{\Delta \cap 
\mathcal{D}_{R}(z)}|\langle \pi(w)\psi_{n}^{\alpha },\pi(z)\psi_n^\alpha\rangle| d\mu ^{+}(w)\leq   \frac{\rho
(\Delta ,R)}{C_{n}^{\alpha }(R)}\text{.}
\end{equation*}
The result thus holds for $p=1$. For $p=\infty$ one trivially has 
\begin{equation*}
\sup_{f\in \mathcal{C}o\, {L}^{\infty }(\mathbb{C}^{+})}\frac{\big\|
W_{\psi _{n}^{\alpha }}f\big|_{\Delta }\big\|_{{L}^{\infty(\C^+) }}}{\Vert W_{\psi _{n}^{\alpha }}f\Vert _{ {L}^{\infty}(\C^+) }}=1.
\end{equation*}
  {Since the coorbit space $\mathcal{C}o\, {L}^{p}(\mathbb{C}^{+})$, $1<p<\infty$, is the complex interpolation space $\big(\mathcal{C}o\, {L}^{\infty}(\mathbb{C}^{+}), \mathcal{C}o\, {L}^{1}(\mathbb{C}^{+})\big)_{1/p}$, see \cite[Theorem~4.7]{fegr89}, the result for $1<p<\infty$ follows from an application of \cite[Theorem~4.1.2]{interpolation}.}  \hfill $\Box $

\begin{corollary}
\label{uncertainty} Let   {$\alpha>1$,} and suppose that $f\in \mathcal{C}o\, {L}^{p}(\mathbb{C}^{+})$, $1\leq p<\infty $, satisfies $\Vert W_{\psi _{n}^{\alpha }}f\Vert
_{L^p(\C^+)}=1$, and that $W_{\psi _{n}^{\alpha }}f$ is $\varepsilon$-concentrated
on $\Delta \subset \mathbb{C}^{+}$, {i.e.}, 
\begin{equation*}
1-\varepsilon \leq \int_{\Delta }|W_{\psi _{n}^{\alpha }}f(z)|^{p}d\mu
^{+}(z)\text{.}
\end{equation*}
Then 
\begin{equation*}
1-\varepsilon \leq \inf_{0<R<1}\left( \frac{\rho (\Delta ,R)}{C_{n}^{\alpha
}(R)}\right) \leq |\Delta |_h.
\end{equation*}
\end{corollary}

\begin{remark}
Let us introduce the Bergman spaces on the unit disk by
\begin{equation}\label{def:berg-disk}
A_{\alpha }^p(\mathbb{D}):=\left\{ f\text{ analytic}:\ \|f\|_{A_\alpha^p}^p:=
\int_{\mathbb{D}}|f(z)|^p(1-|z|^{2})^{\alpha}dz<\infty \right\} .
\end{equation}
We derived   estimates similar to Theorem~\ref{main} for the Bergman 
space $A_{\alpha-2 }^2(\mathbb{D})$ in \cite[Theorem~2]{abspe19}.  
There Seip's double orthogonality result \cite{sei91} was  applied to obtain the concentration estimate.
\end{remark}

\section{Optimal bounds in the analytic setting}\label{sec:opt}


The following result was proved recently for the case $p=2$ in \cite{RamosTilli}. We note here that the case of general $p$ works with exactly the same arguments as Ramos and Tilli used with one extra information: $\mathcal{C}o\,{L}^{p}(\mathbb{C}^{+})$ can be mapped isometrically and surjectively onto a certain Bergman space. 

\begin{theorem}[Ramos-Tilli]
\label{extremal-conjecture} Let   {$1< p<\infty$,} $\Delta \subset \mathbb{C}^{+}$ be a set of
finite hyperbolic measure $A>0$, and $\alpha>1$. Then 
\begin{equation*}
\sup_{|\Delta |_h=A}\sup_{f\in \mathcal{C}o\, {L}^{p}(\mathbb{C}^{+})}%
\frac{\big\| W_{\psi _{0}^{\alpha } }f\big|_{\Delta }\big\|_{L^p(\C^+)}^{p}}{\Vert W_{\psi _{0}^{\alpha }
}f\Vert _{L^p(\C^+)}^{p}}
\end{equation*}%
is attained if and only if $\Delta $ is a pseudohyperbolic disc centered at
some $z^\ast\in \mathbb{C}^+$, up to perturbations of Lebesgue measure zero,  {and $f=\pi(z^\ast)\psi_0^\alpha$}.  In particular,
\begin{equation} \label{eq:ramos-tilli-p}
\int_{\Delta }|W_{\psi _{0}^{\alpha }}f(z)|^{p}d\mu ^{+}(z)\leq  {\left(1-\left(
1+\frac{|\Delta |_h}{4\pi}\right)  ^{1-(\alpha+1) p/2}\right)}\|W_{\psi_0^\alpha}f\|_{L^p(\C^+)}^p\text{.} 
\end{equation}
\end{theorem}

To prove the result for general $p$, it just takes a simple observation. As in \cite{RamosTilli} one can map the space $A_{{(\alpha +1)p/2-2}}^{p}(\mathbb{C}%
^{+})$ to $A_{{(\alpha +1)p/2-2}}^{p}(\mathbb{D})$ (see definitions \eqref{bergman:halfplane} and \eqref{def:berg-disk}) by virtue of the
transformation%
\begin{equation*}
T_{\alpha }F(u)=F\left( i\frac{1+u}{1-u}\right) \left( \frac{1}{1-u}\right)
^{\alpha +1},\quad u\in \mathbb{D}, 
\end{equation*}
which, by \eqref{eq:change-of-variable}, is a multiple of an isometry. Thus,
the composition $T_\alpha\circ s^{-\frac{\alpha +1}{2}}W_{\psi_0^\alpha}:\mathcal{C}o\, {L}^{p}(\mathbb{C}^{+}) \to A_{(\alpha
+1)p/2-2}^{p}(\mathbb{D})$\ 
\begin{equation*}
\mathcal{C}o\, {L}^{p}(\mathbb{C}^{+})  \hspace{0.2cm}\xrightarrow{\hspace{0.2cm}s^{-\frac{\alpha +1}{2}}W_{\psi_0^\alpha }\hspace{0.2cm}}  \hspace{0.2cm}A_{(\alpha +1)p/2-2}^{p}(\mathbb{C}^{+})\hspace{0.2cm}  \xrightarrow{\hspace{0.5cm}T_\alpha\hspace{0.5cm}} \hspace{0.2cm} A_{(\alpha+1)p/2-2}^{p}(\mathbb{D})
\text{,}
\end{equation*}%
is surjective. To see this, we note that by \eqref{eq:n=0} the sequence $\{\psi _{m}^{\alpha }\}_{m\in\N_0}$ (which is complete in $\mathcal{C}%
o\, {L}^{p}(\mathbb{C}^{+})$) is  first mapped to the sequence 
\begin{equation*}
s^{-\frac{\alpha +1}{2}}W_{\psi _{0}^{\alpha }}\psi _{m}^{\alpha }(z)= \sqrt{%
\frac{\Gamma (m+\alpha +1)}{\Gamma (\alpha +1)m!}}\left( \frac{z-i}{z+i}%
\right) ^{m}\left( \frac{2}{1-iz}\right) ^{\alpha +1}, \quad m\in\N_0,\ z=x+is,
\end{equation*}
which is in turn mapped to the sequence of 
monomials.
The completeness of the monomials
in $A_{(\alpha +1)p/2-2}^{p}(\mathbb{D})$,  {$1<p<\infty$, is proven, for example, in \cite[Corollary 4 and Theorem 16]{zhu91}.}
 {Therefore, just as in \cite[Theorem~3.1]{RamosTilli}, the problem of optimal concentration in \eqref{eq:ramos-tilli-p} reduces to optimal concentration on the Bergman space $A^p_{(\alpha+1)p/2-2}(\mathbb{D})$. In fact, the result can be proven for $A^p_{\beta}(\mathbb{D})$, $\beta>0$, without any dependence of $\beta$ on $p$. To do so, one may follow the proof of \cite[Theorem~3.1]{RamosTilli} step by step. We do not repeat the argument here as this would exceed the scope of this paper. There are however only two steps  that need to be adapted. First, the definition of the auxiliary function $u$ (see \cite[p.~8]{RamosTilli}) needs to be changed to $u(z):=|F(z)|^p (1-|z|^2)^{\beta+2}$, $F\in A^p_{\beta}(\mathbb{D})$. Second, the proper analogue of equation (3.9) in \cite{RamosTilli} can be found in \cite[Lemma~3.2]{haakan2000theory}: 
\begin{equation}\label{eq:bergman-book}
u(z)\leq \frac{\beta+1}{\pi}\|F\|^p_{A^p_\beta(\mathbb{D})},\quad z\in\mathbb{D}.\end{equation} Note that we use a different normalization of the area measure on $\mathbb{D}$ than \cite{haakan2000theory}.
The optimal concentration bound is attained if equality in \eqref{eq:bergman-book} holds for some $z^\ast\in \mathbb{D}$. Using the correspondence between wavelet coorbit spaces and Bergman spaces as above (choosing $\beta=(\alpha+1)p/2-2$), we define $$F_{z^\ast}(z )=(1-|z |^2)^{-\frac{\beta+2}{p}}\big\langle \pi(T(z^\ast))\psi_0^{(\beta+1)2/p-1},\pi(T(z ))\psi_0^{(\beta+1)2/p-1}\big\rangle  \in A^p_\beta(\mathbb{D})$$  which satisfies
\begin{align*}
\|u\|_\infty&=\sup_{z\in \mathbb{D}}|F_{z^\ast}(z)|^p(1-|z|^2)^{\beta+2}=\sup_{z\in \mathbb{D}}\big|\big\langle \pi(T(z^\ast))\psi_0^{(\beta+1)2/p-1},\pi(T(z))\psi_0^{(\beta+1)2/p-1}\big\rangle \big|^p
\\&=1=\frac{\beta+1}{\pi}\|F_{z^\ast}\|_{A^p_\beta(\mathbb{D})},
\end{align*}
where we used that $$
\frac{\pi}{\beta+1}=4^{-1}\big\|\big\langle\pi(T(z^\ast))\psi_0^{(\beta+2)2/p-1},\pi(\, \cdot\, )\psi_0^{(\beta+2)2/p-1}\big\rangle\big\|_{L^p(\C^+)}=\|F_{z^\ast}\|_{A^p_\beta(\mathbb{D})}, $$ which follows from left-invariance of the Haar measure, \eqref{eq:change-of-variable}, and \eqref{eq:p-norm-n=0}. This concludes the proof. 
} 

If $W_{\psi _{n}^{\alpha }}f$ is $\varepsilon $
-concentrated on $\Delta \subset \mathbb{C}^{+}$, i.e., 
\begin{equation*}
(1-\varepsilon)\|W_{\psi_0^\alpha}f\|_{L^p(\C^+)}^p \leq \int_{\Delta }|W_{\psi _{0}^{\alpha }}f(z)|^{p}d\mu
^{+}(z)\text{,}
\end{equation*}
then \eqref{eq:ramos-tilli-p} yields 
$$
1-\varepsilon\leq {1-\left(
1+\frac{|\Delta |_h}{4\pi}\right)  ^{1-(\alpha+1) p/2}}.
$$
This leads to the following sharp uncertainty principle for the wavelet
transform.
\begin{theorem}
Let $\alpha>1$, $\varepsilon\in (0,1)$ and $1<
p<\infty $. Suppose that $f\in \mathcal{C}o\, {L}^{p}(\mathbb{C}^{+})$ satisfies $\Vert W_{\psi _{0}^{\alpha }}f\Vert _{L^p(\C^+)}=1$, and that $W_{\psi _{0}^{\alpha }}f$ is $\varepsilon $-concentrated with respect to the $p$-norm in $\Delta \subset 
\mathbb{C}^{+}$. Then 
\begin{equation*}
4\pi\left(\epsilon ^{\frac{2}{2-(\alpha+1) p}}-1\right)\leq |\Delta |_h.
\end{equation*}
\end{theorem}

\begin{remark}
Extending these optimal results to the general family of wavelets $\psi _{n}^{\alpha
}$ remains an open problem (as well as the corresponding optimal bound for
general Hermite functions in \cite{abspe18-sieve}). The methods of 
\cite{Inventiones} and \cite[Theorem~3.1]{RamosTilli} are extremely
dependent on the analytic structure of the case $n=0$, and for $n>0$ the transforms are not related to analytic functions. However, since a derivative of an analytic function is still analytic, the methods of 
\cite{Inventiones} and \cite[Theorem~3.1]{RamosTilli}  can be adapted to obtain sharp inequalities for weighted $L^p$-norms of derivatives of holomorphic functions.   In \cite%
{kalaj2022contraction}, combined with sharp pointwise inequalities for derivatives in Fock spaces (involving Laguerre functions with complex argument), this has been used to obtain a contraction inequality for derivatives
in Fock spaces. It is reasonable to expect a related inequality for derivatives in Bergman spaces. This requires sharp pointwise inequalities for derivatives in Bergman spaces, which possibly involve Jacobi polynomials with complex argument.
\end{remark}

\subsection{Lieb's uncertainty principle}

The following result is due to Kulikov \cite[Theorem~1.2]{kuli22}. 
Note that the normalization of the measure on $ \mathbb{D}$   is chosen such that $\Vert 1\Vert^p_{A_{\alpha
 }^{p}(\mathbb{D})}= {\frac{\pi}{\alpha+1}}$ (as defined in \eqref{def:berg-disk}).

\begin{theorem}[Kulikov]
\label{thm:kulikov} Let $G:[0,\infty )\rightarrow \mathbb{R}^{+}$ be a
convex function. Then for every $f\in A_{\alpha-2 }^{p}(\mathbb{D})$ with $%
\Vert f\Vert^p_{A_{\alpha-2 }^{p}(\mathbb{D})}= {\frac{\pi}{\alpha+1}}$ it holds 
\begin{equation*}
\int_{\mathbb{D}}G\left( |f(z)|^{p}(1-|z|^{2})^{\alpha }\right)
(1-|z|^{2})^{-2}dz\leq \int_{\mathbb{D}}G\left( (1-|z|^{2})^{\alpha }\right)
(1-|z|^{2})^{-2}dz,
\end{equation*}%
and equality  is attained for $f\equiv 1$.
\end{theorem}

The next theorem is simply a restatement of Kulikov's theorem for a particular choice of $G$, but we
would like to highlight this case since this is the first instance of a counterpart
of Lieb's inequality for the short-time Fourier transform. From the
experience in time-frequency analysis, this inequality should have many
applications (as noticed in the comments after Corollary 1.2 in \cite{nicola2022norm}, Lieb's inequality is equivalent to sharp norms for localization operators; for applications in signal analysis see \cite{baraniuk2001measuring} and \cite{ricaud2014survey}). It will also be essential in the derivation of the local Lieb
formula in the next section. 

\begin{theorem}[Lieb's inequality for analytic wavelets]
\label{thm:Lieb-UP}Let $\alpha >1$,   
  $p\geq 2$, and $f\in H^2(\C^+)$. Then
\begin{equation}
{\Vert W_{\psi _{0}^{\alpha }}f\Vert _{ {L}^{p}(\mathbb{C}^{+})}^{p}}%
\leq \frac{8\pi  }{ {  (\alpha  +1 )p-2}}\Vert f\Vert_{H^2(\C^+) }^{p},
\end{equation}
with the inequality being sharp.
\end{theorem}

\noindent\textbf{Proof:\ } If we set $\alpha =  {\beta +1}$ and define 
\begin{equation*}
F(z)=s^{-(\beta /2+1)}\left( \Vert W_{\psi _{0}^{ {\beta +1}}}f\Vert _{%
{L}^{2}(\mathbb{C}^{+})}\right) ^{-1}\ W_{\psi _{0}^{ {\beta +1}}}f(z),
\end{equation*}%
then $F$ is analytic and $\Vert F\Vert _{{A}^{2}_\beta(\mathbb{C}%
^{+})}=1$. Moreover, taking $G(t)=t^{p/2}$ (which is convex for $p\geq 2$%
) in Theorem~\ref{thm:kulikov} and applying \eqref{eq:change-of-variable}
shows 
\begin{align*}
& \frac{\Vert W_{\psi _{0}^{{\beta +1}}}f\Vert _{ {L}^{p}(\mathbb{C}%
^{+})}^{p}}{\Vert W_{\psi _{0}^{ {\beta +1}}}f\Vert _{{L}^{2}(%
\mathbb{C}^{+})}^{p}}=\int_{\mathbb{C}^{+}}|F(z)|^{p}s^{\beta p/2+p-2}d\mu_{\C^+}(z) \\
& \hspace{1cm}=4\left( \frac{\beta +1}{4\pi }\right) ^{\frac{p}{2}}\int_{%
\mathbb{D}}\left( \frac{4\pi }{\beta +1}|F(T(u))|^{2}|1-u|^{-(2\beta
+4)}(1-|u|^{2})^{\beta +2}\right) ^{\frac{p}{2}}(1-|u|^{2})^{-2}du \\
& \hspace{1cm}\leq 4\left( \frac{\beta +1}{4\pi }\right) ^{\frac{p}{2}}\int_{%
\mathbb{D}}(1-|u|^{2})^{\beta p/2+p-2}du=4\pi \left( \frac{\beta +1}{4\pi }%
\right) ^{\frac{p}{2}}\int_{0}^{1}r^{\beta p/2+p-2}dr \\
& \hspace{1cm}=\frac{4\pi }{p(\beta /2+1)-1}\left( \frac{\beta +1}{4\pi }%
\right) ^{\frac{p}{2}}.
\end{align*}%
Note that we were allowed to apply Theorem~\ref{thm:kulikov} as 
\begin{align*}
\left\Vert  {\frac{4\pi }{\beta +1}}F(T(u))^{2}(1-u)^{-(2\beta +4)}\right\Vert
_{A_{\beta }^{1}(\mathbb{D})}& = \frac{4\pi }{\beta +1}\int_{\mathbb{D}}|F(T(u))|^{2}|1-u|^{-(2%
\beta +4)}(1-|u|^{2})^{\beta }du \\
& =  {\frac{\pi }{\beta +1}}\int_{\mathbb{C}^{+}}|F(z)|^{2}s^{\beta }dz=  {\frac{\pi }{\beta +1}}.
\end{align*}%
Substituting back $\alpha =  {\beta +1}$ and using $\Vert W_{\psi }f\Vert _{{L}^{2}(\mathbb{C}^{+})}^{2}=\|f\|^2_{H^2(\C^+)}\|\widehat{\psi}\|_{L^2(\R^+,t^{-1})}^2$ yields 
\begin{align*}
\Vert W_{\psi _{0}^{\alpha }}f\Vert _{{L}^{p}(\mathbb{C}^{+} )}^{p}&
\leq \frac{4\pi  }{ {p (\alpha  +1 )/2-1}}\left( { \frac{\alpha}{4\pi }}\right)
^{p/2}\Vert f\Vert_{H^2(\C^+)} ^{p}\big\| \widehat{\psi _{0}^{\alpha }}\big\| _{L^{2}(%
\mathbb{R}^{+},t^{-1})}^{p} \\
& =\frac{4\pi  }{ {p (\alpha  +1 )/2-1}} \Vert f\Vert ^{p}_{H^2(\C^+)},
\end{align*}%
where we used that $\big\| \widehat{\psi _{0}^{\alpha }}\big\| _{L^{2}(\mathbb{%
R}^{+},t^{-1})}^{2}=4\pi/\alpha $. \hfill $
\Box $\newline

\subsection{Local Lieb's uncertainty principle}
We conclude our considerations on sharp concentration inequalities with the following immediate consequence of the results of the last two subsections. This is the wavelet analogue of Theorem 5.2 in \cite{Inventiones}, but we provide a more direct proof.
\begin{theorem}[Local Lieb's inequality for analytic wavelets]
Let $\Delta \subset \mathbb{C}^{+}$ be a set of finite measure, $\alpha>1$, and $2\leq p<\infty$. For $f\in H^2(\C^+)$, it holds
\begin{equation*}
\frac{\big\| W_{\psi_0^\alpha }f\big|_{\Delta }\big\| _{L^p(\C^+) }^{p}}{\Vert f\Vert _{H^2(\C^+)}^{p}}\leq\frac{8\pi  }{{ (\alpha  +1 )p-2}} {\left(1-\left(
1+\frac{|\Delta |_h}{4\pi}\right)  ^{1-(\alpha+1) p/2}\right)}.
\end{equation*}
\end{theorem}
\noindent\proof Apply first \eqref{eq:ramos-tilli-p} and then  Theorem~\ref{thm:Lieb-UP}.\pbox


\subsection{Discussion of general Lieb inequalities}

The currently available methods do not seem to provide a Lieb inequality
with sharp constants for a general window $\psi $. The problem can be stated
as providing the best constant in the inequality%
\begin{equation}
\Vert W_{\psi }f\Vert _{ {L}^{p}(\mathbb{C}^{+})}^{p}\leq C(\psi
,p)\Vert f\Vert ^{p}_{H^2(\C^+)}\text{.}  \label{Lieb}
\end{equation}%
This seems to be a non-trivial question.  Kulikov's methods are dependent on
analytic functions, which, by the results of \cite{AnalyticWavelet}, restricts
the window to a weight times $\psi _{0}^{\alpha }$. On the other side, Lieb's 
methods depend on the optimal constants of Young inequality, which is not possible to conveniently use in this context.
An intermediate result
would be to obtain \eqref{Lieb} for the family $\psi _{n}^{\alpha }$, which
still keeps some properties of the analytic case.\ We will obtain an
inequality for $\Vert W_{\psi }f\Vert _{{L}^{p}(\mathbb{C}^{+})}^{p}$
using the information from Theorem \ref{thm:Lieb-UP}, but our estimates are
quite rough, and serve only to give an idea of a significant obstruction we
face, namely the non-unimodularity of the affine group, which does not allow
to apply Young's inequality in a conventional form (there are versions of Young's inequality for the affine group, see Section 4 in \cite{klein1978sharp}, but they do not seem to apply to the problem at hand).  {For locally compact Abelian groups, Lieb's uncertainty principle is well-established, and its optimizers were recently characterized in \cite{nicola-lca}.}

For any non-unimodular locally compact group $\mathcal{G}$ with left-invarian Haar measure $\mu_{\mathcal{G}}$, we show the
following version of Young's convolution inequality. First, the convolution
of two functions $F,G:\mathcal{G}\rightarrow \mathbb{C}$ is defined as 
\begin{equation*}
(F\ast G)(z)=\int_{\mathcal{G}}F(w)G(w^{-1}\cdot z)d\mu_{\mathcal{G}}(w).
\end{equation*}
In addition, we define reflection operator $\mathcal{R}F(z)=F(z^{-1})$.

\begin{lemma}
\label{lem:young} Let $\mathcal{G}$ be a locally compact group with left
Haar measure $\mu_{\mathcal{G}}$, and $1\leq p,q,r\leq \infty$ be such that $1+\frac{1}{r}=%
\frac{1}{p}+\frac{1}{q}$. If $F\in L^p(\mathcal{G})$, and $G, \mathcal{R}%
G\in L^q(\mathcal{G})$, then 
\begin{equation*}
\|F\ast G\|_{L^r(\mathcal{G})}\leq \|F\|_{L^p(\mathcal{G})}\cdot \max\{\|G\|_{L^q(\mathcal{G})},\|\mathcal{R}G\|_{L^q(\mathcal{G})}\}. 
\end{equation*}
\end{lemma}

\noindent \textbf{Proof:\ } We adopt the classical proof of Young's inequality for
unimodular groups via the Riesz-Thorin theorem applied to the operator $T_G
F=F\ast G$. Let us assume that $G, \mathcal{R}G\in L^q(\mathcal{G})$. First,
if $r=\infty$ we note that by H\"older's inequality 
\begin{equation*}
\|T_GF\|_{L^\infty(\mathcal{G})} \leq \sup_{z\in \mathcal{G}} \int_{\mathcal{G}}
|F(w)\mathcal{R}G(z^{-1}\cdot w)|d\mu_{\mathcal{G}}(w)\leq \|F\|_{L^p(\mathcal{G})}\|\mathcal{R}G\|_{L^q(\mathcal{G})},\quad \frac{1}{p%
}+\frac{1}{q}=1. 
\end{equation*}
On the other hand, for $r=q$, we have   that $\|T_G F\|_{L^q(\mathcal{G})}\leq \|F\|_{L^1(\mathcal{G})}\|  G\|_{L^q(\mathcal{G})}$ by   \cite[(20.14)]{hewittross}. An application of the 
Riesz-Thorin theorem then completes the proof. \hfill$\Box$\newline

Together with Theorem \ref{thm:Lieb-UP} this leads to the following.

\begin{proposition}
Let $f\in H^2(\C^2)$, $\psi\in \mathcal{C}o L^1(\C^+)$ be an admissible wavelet. For every $\alpha >1,$ and every $2\leq p<\infty$, it holds 
\begin{align*}
\frac{\Vert W_{\psi }f\Vert _{{L}^{p}(\mathbb{C}^{+})}^{p}}{\|f\|^p_{H^2(\C^+)}}\leq
 {\frac{8\pi  }{ (\alpha  +1 )p-2}\left( \frac{\alpha }{4\pi}\right)
^{p}}\max\left\{ {\Vert W_{\psi }{\psi _{0}^{\alpha }}\Vert
_{{L}^{1}(\mathbb{C}^{+})}^{p}},\Vert W_{{{\psi _{0}^{\alpha }}}%
}\psi \Vert _{{L}^{1}(\mathbb{C}^{+})}^{p}\right\}. 
\end{align*}
\end{proposition}

\noindent\proof  Let us  apply consecutively \eqref{ortogonalityrelations},  Lemma~\ref{lem:young},  and  Theorem~\ref{thm:Lieb-UP}  to show 
\begin{align*}
\Vert W_{\psi}f \Vert _{{L}^{p}(\mathbb{C}^{+})}^{p} &=\Vert 
\widehat{\psi _{0}^{\alpha }}\Vert _{L^{2}(\mathbb{R}^{+},t^{-1})}^{-2p}
\int_{\mathbb{C}^{+}}\left\vert \int_{\mathbb{C}^{+}}W_{\psi _{0}^{\alpha
}}f(w)\langle {{\psi _{0}^{\alpha }}},\pi (w^{-1}\cdot z)\psi\rangle d\mu^+(w)\right\vert ^{p}d\mu^+(z) 
\\
& \leq  \Vert \widehat{\psi _{0}^{\alpha }}\Vert _{L^{2}(\mathbb{R}%
^{+},t^{-1})}^{-2p}\Vert W_{\psi _{0}^{\alpha }}f\Vert _{ {L}^{p}(%
\mathbb{C}^{+})}^{p}\max \left\{ \Vert W_{\psi}{{\psi _{0}^{\alpha }}}\Vert _{%
{L}^{1}(\mathbb{C}^{+} )}^{p},\Vert W_{{{\psi _{0}^{\alpha }}}}\psi
\Vert _{ {L}^{1}(\mathbb{C}^{+} )}^{p}\right\} 
\\
& \leq  {\frac{4\pi  }{p (\alpha  +1 )/2-1}\hspace{-2pt}\left( \frac{\alpha }{4\pi}\right)
^{p}}\Vert f\Vert ^{p}_{H^2(\C^+)}\max \hspace{-1pt}\left\{ {\Vert W_{\psi }{\psi _{0}^{\alpha }}\Vert
_{{L}^{1}(\mathbb{C}^{+} )}^{p}},\Vert W_{{{\psi _{0}^{\alpha }}}%
}\psi \Vert _{{L}^{1}(\mathbb{C}^{+} )}^{p}\right\} ,
\end{align*}%
where we used that $\mathcal{R}W_{\psi}f=\overline{W_{f}\psi}$. \pbox

\section{Appendix}

\subsection{Calculation of Reproducing Kernels and Basis Functions}

The following integral formula \cite[p. 809, 7.414 (4)]{grary07} enables us
to explicitly calculate the wavelet transform $W_{\psi _{n}^{\alpha }}$ of
the basis elements $\psi _{m}^{\alpha }$ of $H^{2}(\mathbb{C}^{+})$ 
\begin{align}
\int_{0}^{\infty }e^{-bt}t^{\alpha }L_{n}^{\alpha }(\lambda t)L_{m}^{\alpha
}(\mu t)dt =&\frac{\Gamma (m+n+\alpha +1)}{m!n!}\frac{(b-\lambda )^{n}(b-\mu
)^{m}}{b^{m+n+\alpha +1}}  \notag  \label{grad-equ} \\
& \times F\left( -m,-n;-m-n-\alpha ;\frac{b(b-\lambda -\mu )}{(b-\mu
)(b-\lambda )}\right) ,
\end{align}%
where $\mbox{Re}(\alpha )>-1$, $\mbox{Re}(b)>0$, and $F={}_{2}F_{1}$ denotes
the \emph{Gauss hypergeometric} function defined by 
\begin{equation}
F(a,b;c;z)=\sum_{k=0}^{\infty }\frac{(a)_{k}(b)_{k}}{(c)_{k}k!}z^{k},
\label{eq:hypergeom}
\end{equation}%
and $(x)_{k}$ denotes Pochhammer's rising factorial. Moreover, we need the
following explicit representation of the \emph{Jacobi polynomials} \cite[p.
442, (18.5.7)]{NIST10} 
\begin{equation}
P_{n}^{(\alpha ,\beta )}(z)=\frac{\Gamma (n+\alpha +1)}{n!\Gamma (n+\alpha
+\beta +1)}\sum_{k=0}^{n}\binom{n}{k}\frac{\Gamma (n+k+\alpha +\beta +1)}{
\Gamma (k+\alpha +1)}\left( \frac{z-1}{2}\right) ^{k}.  \label{eq:jacobi}
\end{equation}
This formula implies the following connection between the hypergeometric
function and the Jacobi polynomials.

\begin{lemma}
\label{lem:jacobi-hypergeo} For $n,m\in\mathbb{N}_0$, it holds 
\begin{align*}  \label{eq:jacobi-hypergeo}
F(-n,&-m;-n-m-\alpha;z) \\
&= \dfrac{\min(n,m)!\Gamma(\max(n,m)+\alpha+1)}{\Gamma(n+m+\alpha+1)}
(-z)^{\min(n,m)}P_{\min(n,m)}^{(|n-m|,\alpha)}\left(1-\frac{2}{z}\right).
\end{align*}
\end{lemma}

\noindent\textbf{Proof:\ } As $F(-n,-m;-n-m-\alpha;z)=F(-m,-n;-n-m-\alpha;z)$, we may
assume that $n\geq m$. The definition of $F$ and   $
(-x)_k=(-1)^k(x-k+1)_k $ lead us to
\begin{align*}
F(-n,-m;&-n-m-\alpha;z)=\sum_{k=0}^\infty \frac{(-n)_k(-m)_k}{
(-n-m-\alpha)_k k!}z^k \\
&=\sum_{k=0}^m \frac{(n-k+1)_k(m-k+1)_k}{(n+m+\alpha-k+1)_k k!}(-z)^k \\
&=\sum_{k=0}^m \frac{\Gamma(n+1)\Gamma(m+1)\Gamma(n+m+\alpha-k+1)}{
\Gamma(n-k+1)\Gamma(m-k+1)\Gamma(n+m+\alpha+1) k!}(-z)^k \\
&=\frac{n!}{\Gamma(n+m+\alpha+1)}\sum_{k=0}^m \binom{m}{k}\frac{
\Gamma(n+m+\alpha-k+1)}{\Gamma(n-k+1)}(-z)^k \\
&=\frac{n!}{\Gamma(n+m+\alpha+1)}\sum_{k=0}^m \binom{m}{k}\frac{
\Gamma(n+\alpha+k+1)}{\Gamma(n-m+k+1)}(-z)^{m-k} \\
&=\frac{\Gamma(n+1)}{\Gamma(n+m+\alpha+1)}(-z)^m\sum_{k=0}^m \binom{m}{k}
\frac{\Gamma(n+\alpha+k+1)}{\Gamma(n-m+k+1)}\left(-\frac{1}{z}\right)^{k} \\
&=\frac{m!\Gamma(n+\alpha+1)}{\Gamma(n+m+\alpha+1)}(-z)^m
P_m^{(n-m,\alpha)}\left(1-\frac{2}{z}\right).
\end{align*}
\hfill$\Box$\newline

\noindent \textbf{Proof of Proposition \ref{lem:explicit-wavelet-trafo}:} Let $z=x+is$ and $w=x^{\prime }+is^{\prime }$. By definition of $
\psi _{n}^{\alpha }$, one has 
\begin{align*}
\frac{1}{2^{\alpha +2}\pi}& \sqrt{\frac{\Gamma (n+\alpha +1)\Gamma (m+\alpha +1)
}{n!m!}}\langle \pi (w)\psi _{m}^{\alpha },\pi (z)\psi _{n}^{\alpha }\rangle 
\\
& \hspace{2.5cm}=  {\frac{1}{2\pi}}\sqrt{ss^{\prime }}\int_{0}^{\infty }(s^{\prime }t)^{\frac{\alpha }{2}}e^{-s^{\prime }t-ix^{\prime }t}L_{m}^{\alpha }(2s^{\prime
}t)(st)^{\frac{\alpha }{2}}e^{-st+ixt}L_{n}^{\alpha }(2st){dt} \\
& \hspace{2.5cm}= {\frac{1}{2\pi}}(ss^{\prime })^{\frac{\alpha +1}{2}}\int_{0}^{\infty
}t^{\alpha }e^{it(z-\overline{w})} L_{n}^{\alpha }(2st)L_{m}^{\alpha }(2s^{\prime
}t)dt=\circledast .
\end{align*}
Hence, setting $b=-i(z-\overline{w})$, $\lambda =2s$, and $\mu =2s'$ in
equation \eqref{grad-equ} yields $b-\lambda =-i(\overline{z}-\overline{w})$, $b-\mu =-i(z-w)$, $b-\lambda -\mu =-i(\overline{z}-w)$, and consequently 
\begin{align*}
\circledast =& {\frac{1}{2\pi}} (ss^{\prime })^{\frac{\alpha +1}{2}}\frac{\Gamma (m+n+\alpha
+1)}{m!n!}\frac{(z-w)^{m}(\overline{z}-\overline{w})^{n}}{(z-\overline{w}
)^{m+n}}\big(-i(z-\overline{w})\big)^{-\alpha -1} \\
& \times F\left( -m,-n;-m-n-\alpha ;\left\vert \frac{z-\overline{w}}{z-w}
\right\vert ^{2}\right) .
\end{align*}
Setting $n=m$ gives \eqref{eq:kernel-expl} by Lemma~\ref{lem:jacobi-hypergeo}. Taking $w=i$ and applying Lemma~\ref{lem:jacobi-hypergeo} again yields \eqref{eq:explicit-wavelet}. {Finally,   \eqref{C_psi_n^alpha} follows from the basic identity $L_n^{\alpha+1}=\sum_{k=0}^nL_k^\alpha$ and the orthogonality relation for generalized Laguerre polynomials (see, e.g., \cite{NIST10}) \
\begin{align*}
    C_{\psi_n^\alpha}&=\|\widehat\psi_n^\alpha\|^2_{L^2(\R^+,t^{-1})}=\frac{2^{\alpha+2}\pi n!}{\Gamma(n+\alpha+1)}\int_0^\infty t^{\alpha-1}e^{-2t}L_n^\alpha(2t)^2dt
    \\
    &=\frac{4\pi n!}{\Gamma(n+\alpha+1)} \int_0^\infty t^{\alpha-1}e^{-t}\left(\sum_{k=0}^nL_k^{\alpha-1}(t)\right)^2dt
        \\
    &=\frac{4\pi n!}{\Gamma(n+\alpha+1)} \sum_{k=0}^n\int_0^\infty t^{\alpha-1}e^{-t} L_k^{\alpha-1}(t)^2dt=\frac{4\pi n!}{\Gamma(n+\alpha+1)}\sum_{k=0}^n\frac{\Gamma(k+\alpha)}{k!}=\frac{4\pi}{\alpha},
\end{align*}
where the last identity can be easily shown by an inductive argument.
}  
\hfill $\Box 
$\newline

\subsection{Integrability}\label{app:int}

It remains for us to show that the mother wavelets $\psi _{n}^{\alpha }$
are admissible for the integrable representation $\pi $ of the $ax+b$ group for appropriate choices of $\alpha$. 



\begin{proposition}
\label{prop:integrable} Let $n\in \N_0$. If  $\alpha >1$, then 
\begin{equation}
\int_{\mathbb{C}^{+}}\big|\langle \psi _{n}^{\alpha },\pi (z)\psi
_{n}^{\alpha }\rangle \big|d\mu^+(z)<\infty .  \label{eq:integrable}
\end{equation}
{Moreover,  if $(\alpha+1)p-2>0$ then
\begin{equation}\label{eq:p-norm-n=0}
\int_{\mathbb{C}^{+}}\big|W_{\psi _{0}^{\alpha }}\psi _{0}^{\alpha }(z)\big|^p
\frac{d\mu_{\C^+}(z)}{\emph{Im}(z)^{2}}=\frac{8\pi}{(\alpha+1)p-2}.
\end{equation}
}
\end{proposition}

\noindent\textbf{Proof:\ } By  \eqref{eq:aux-moebius}, we have 
\begin{equation*}
\big|W_{\psi _{n}^{\alpha }}\psi _{n}^{\alpha }(T(u))\big|=|u|^{2n}\left( 
\frac{1-|u|^{2}}{4}\right) ^{\frac{\alpha +1}{2}}|Z_{n,n}^\alpha(|u|^{2})|.
\end{equation*}%
Consequently, by \eqref{eq:change-of-variable} and \eqref{eq:def-Phi_nm}
\begin{align*}
\int_{\mathbb{C}^{+}}\big|W_{\psi _{n}^{\alpha }}\psi _{n}^{\alpha }(z)\big|
\frac{d\mu_{\C^+}(z)}{\text{Im}(z)^{2}}& = \int_{\mathbb{D}}\big|W_{\psi _{n}^{\alpha }}\psi
_{n}^{\alpha }(T(u))\big|\frac{4du}{(1-|u|^{2})^{2}} \\
& =4^{-\frac{\alpha -1}{2}}\int_{\mathbb{D}}|u|^{2n}|Z
_{n,n}^\alpha(|u|^{2})|(1-|u|^{2})^{\frac{\alpha -3}{2}}du \\
& \lesssim \int_{\mathbb{D}}(1-|u|^{2})^{\frac{\alpha -3}{2}}du=\pi
\int_{0}^{1}r^{\frac{\alpha -3}{2}}dr.
\end{align*}
The last integral is finite if and only if $\alpha >1$. {For $n=0$ and $p\geq 1$ we argue similarly to obtain
\begin{align*}
\int_{\mathbb{C}^{+}}\big|W_{\psi _{0}^{\alpha }}\psi _{0}^{\alpha }(z)\big|^p
\frac{d\mu_{\C^+}(z)}{\text{Im}(z)^{2}}& = \int_{\mathbb{D}}\big|W_{\psi _{0}^{\alpha }}\psi
_{0}^{\alpha }(T(u))\big|^p\frac{4du}{(1-|u|^{2})^{2}} \\
& =4 \int_{\mathbb{D}}   (1-|u|^{2})^{\frac{(\alpha +1)p}{2}-2}du \\
& =4 \pi\int_{0}^1(1-t)^{\frac{(\alpha +1)p}{2}-2}dt
\\
&=4 \pi
\frac{1}{\frac{(\alpha+1)p}{2}-1}=\frac{8\pi}{(\alpha+1)p-2}.
\end{align*}}
\hfill $\Box $
\newline

\section*{Acknowledgements}
We would like to thank the anonymous reviewers. Their valuable input helped  to substantially improve this article.
The authors ackowledge the support of the Austrian Science Fund (FWF)
through the projects 10.55776/P31225 (L.D. Abreu) and 10.55776/J4254, and 10.55776/Y1199 (M. Speckbacher).


\bibliographystyle{plain}
\bibliography{paperbib}

\end{document}